%% file: 2005-14.tex
\documentclass{gtart_h}  

\input gtoutput

\lognumber{550}
\received{30 December 2004}
\volumenumber{9}\papernumber{14}\volumeyear{2005}
\pagenumbers{483}{491}   
\published{7 April 2005}
\accepted{25 March 2005}
\proposed{Yasha Eliashberg}
\seconded{Leonid Polterovich, Simon Donaldson}

\usepackage{amssymb,amsmath}

\newtheorem{lemma}{Lemma}
\newtheorem{theorem}{Theorem}
\newtheorem{corollary}{Corollary}
\newtheorem{conjecture}{Conjecture}
\theoremstyle{definition}

\newtheorem{remark}{Remark}

\newcommand{\Z}{{\mathbb Z}}
\newcommand{\R}{{\mathbb R}}

\newcommand{\PP}{{\mathbb P}}

\newcommand{\Proj}{{\mathbf{Proj}}}

\let\proofend\endproof

\newcommand{\Pic}{{\operatorname{Pic}}}

\newcommand{\GW}{{\operatorname{GW}}}

\newcommand{\Sig}{{\Sigma}}

\begin{document}
\title[Asymptotics of Gromov--Witten invariants]{Logarithmic asymptotics of the
genus zero\\Gromov--Witten invariants of the blown up plane}

\asciititle{Logarithmic asymptotics of the
genus zero Gromov-Witten invariants of the blown up plane}

\author{Ilia Itenberg\\Viatcheslav Kharlamov\\Eugenii Shustin}
\shortauthors{Ilia Itenberg, Viatcheslav Kharlamov and Eugenii Shustin}
\asciiauthors{Ilia Itenberg, Viatcheslav Kharlamov and Eugenii Shustin}

\addresses{{\rm(II and VK)} Universit\'e Louis Pasteur et IRMA,
7, rue Ren\'e Descartes,\\67084 Strasbourg Cedex, France
\\{\rm and}\\
{\rm(ES)} School of Mathematical Sciences,
Tel Aviv University\\Ramat Aviv, 69978 Tel Aviv, Israel\\
\smallskip\\{\rm Email:}\qua{\tt\mailto{itenberg@math.u-strasbg.fr},
\mailto{kharlam@math.u-strasbg.fr}\\
\mailto{shustin@post.tau.ac.il}}}

\asciiaddress{(II and VK) Universite Louis Pasteur et IRMA,
7, rue Rene Descartes\\67084 Strasbourg Cedex, France\\and\\(ES) 
School of Mathematical Sciences,
Tel Aviv University\\Ramat Aviv, 69978 Tel Aviv, Israel}

\asciiemail{itenberg@math.u-strasbg.fr, kharlam@math.u-strasbg.fr,
shustin@post.tau.ac.il}

\begin{abstract}
We study the growth of the genus zero Gromov--Witten invariants $GW_{nD}$
of the projective plane $P^2_k$ blown up at $k$ points
(where $D$ is a class in the second homology group of $P^2_k$).
We prove that, under some natural restrictions on $D$,
the sequence $\log GW_{nD}$ is equivalent to $\lambda n\log n$,
where $\lambda = D \cdot c_1(P^2_k)$.
\end{abstract}
\asciiabstract{%
We study the growth of the genus zero Gromov-Witten invariants GW_{nD}
of the projective plane P^2_k blown up at k points
(where D is a class in the second homology group of P^2_k).
We prove that, under some natural restrictions on D,
the sequence log GW_{nD} is equivalent to lambda n log n,
where lambda = D.c_1(P^2_k).}

\asciikeywords{Gromov-Witten invariants, rational and ruled algebraic
surfaces, rational and ruled symplectic 4-manifolds, tropical
enumerative geometry} 

\maketitle

\section{Introduction}\label{intro}

In this note we treat the asymptotic behavior of
the genus zero
Gromov--Witten invariants
on $4$--dimensional symplectic manifolds.
In this setting such an invariant can be seen as a count of connected
rational $J$--holomorphic curves in a given homology class
under a choice of a generic tamed almost complex structure,
see, for example, \cite{MS}.

Our interest to logarithmic asymptotics is motivated by a comparison
of Grom\-ov--Witten invariants with their real analogs introduced by
J-Y~Welschinger (see \cite{IKS1} and the conjecture in Section
\ref{last}) and by a relation of the logarithmic asymptotics with the
convergency properties of the Gromov--Witten potential (see
\cite{IKS1}).

As is known, already the existence of homology classes with a
nontrivial invariant which are distinct from the homology classes
of $(-1)$--curves is a very restrictive condition. It implies that
the $4$--dimensional symplectic manifold in question is a blow-up
of a rational or ruled manifold (precise statements, details, and
references can be found in~\cite{MS}, Section 9.4).
We exclude irrational ruled
manifolds (that is, symplectic $S^2$--bundles over Riemann surfaces of
genus $g>0$) since they have only one homology class with a
nontrivial invariant, the class represented by the fiber.
Furthermore, since the Gromov--Witten invariants are preserved
under variations of the symplectic structure, for the study of
their asymptotic properties it is sufficient to consider the
product of complex projective lines, $\PP^1\times\PP^1$, and the
blow-ups of the complex projective plane.

Let us denote by $\PP^2_k$ the complex projective plane $\PP^2$
blown up at $k$ points. Pick a homological class~$D$ in
$H_2(\PP^2_k; \Z)$ such that the Gromov--Witten invariant
$GW_D(\PP^2_k)$ is non-zero, and either $D \cdot c_1(\PP^2_k) >
2$, or $D \cdot c_1(\PP^2_k) = 2$ and
$D^2 > 0$.
Under the
above hypotheses on~$D$, the Gromov--Witten invariants of $nD, n\ge
1,$ are enumerative, that is, the invariant $\GW_{nD} (\PP^2_k),
n\ge 1$, is equal to the number $N_{nD} (\PP^2_k)$ of immersed
irreducible rational curves passing through $nD\cdot c_1(\PP^2_k)
-1$ given generic points in $\PP^2_k$ under the additional
assumption that the blown up points are also generic, see
\cite{GP}.

In the case of $\PP^2_0=\PP^2$
the Kontsevich recursive formula
for $N_{nL} (\PP^2)$
\cite{KM}
($L$ being a line in $\PP^2$) allows one
to get
successive values of these invariants
and
to find their
asymptotics.
In particular, one has
$\log N_{nL}(\PP^2) = 3n\log n+O(n)$
as
$n\to +\infty$
(see \cite{Itz}).
There exist recursive formulas for the Del Pezzo surfaces, see
\cite{KM}, and for $\PP^2_k$ with any $k$,
see \cite{GP}.
However, these formulas
are not easy to
analyze
specially for large $k$.
In \cite{IKS1} working with the corresponding counts of real curves
we observed,
by means of Mikhalkin's
theorem \cite{Mi0,Mi}
(see
also~\cite{Sh0}) on the enumeration of nodal curves on toric
surfaces via lattice paths in convex lattice polygons, that the
relation
$$
\log N_{nD}(\Sig)=\lambda \,n\log n + O(n),\quad \lambda =D\cdot
c_1(\Sig),
$$
holds
for any ample divisor $D$ on a toric Del Pezzo
surface $\Sig$, in particular, on the plane with blown up one,
two, or three points, and on $\PP^1\times \PP^1$.

The present note is devoted to a proof of
the following theorem.

\begin{theorem}\label{nt1}
Let $\PP^2_k$ be the plane blown up
at $k\ge 1$
points, and
$D\in H_2(\PP^2_k;\Z)$
a homology class
such that $GW_D(\PP^2_k)\ne 0$
and either $D \cdot c_1(\PP^2_k) > 2$, or
$D \cdot c_1(\PP^2_k) = 2$ and $D^2 > 0$.
Then
\begin{equation}
\log GW_{nD}(\PP^2_k)= \lambda\, n\log n + O(n),\quad \lambda=D\cdot
c_1(\PP^2_k)\ .
\label{ne1}\end{equation}\end{theorem}

As a consequence we get the following enumerative statement.

\begin{corollary}\label{cnt1}
Let $\PP^2_k$ be the plane blown up at $k\ge 1$ generic points,
and $D\in H_2(\PP^2_k;\Z)$ is as in
{\rm Theorem \ref{nt1}}. Then
\begin{equation}\log N_{nD}(\PP^2_k)= \lambda\, n\log n + O(n),\quad
\lambda=D\cdot c_1(\PP^2_k)\ .\label{ne10}\end{equation}
Furthermore, if $k\le 9$, then {\rm(\ref{ne10})} holds for any ample
divisor $D$ on $\PP^2_k$.
\end{corollary}

Let us notice that the hypotheses of Theorem~\ref{nt1} are in a
sense optimal. For example, $GW_{nD}=0$ if $n\ge 2$ and $D$ is an
embedded curve with $D\cdot c_1\le 2$.

\section{Rational curves on
rational geometrically ruled surfaces}

Here we prove
two auxiliary statements.

\begin{lemma}\label{nl1}
Let $\Sigma_s$,
$s>0$, be a rational geometrically ruled
surface with
the exceptional section $E$,
$E^2=-s$, and a fibre
$F$. Then,
\begin{equation}
\log N_{n(sF+E)} (\Sigma_s) \ge (s + 2)n\log n+O(n)\
.\label{ne2}\end{equation}
\end{lemma}

{\bf Proof}\qua We follow the ideas of the proof of Lemma 5 in
\cite{IKS1}.

First, we observe that the case
$s=1$
corresponds
to curves on $\Sig=\PP^2_1$
disjoint from $E$, and since
$s + 2 = (sF+E)\cdot
c_1(\Sig)$, this case is
settled in
Theorem 3 of
\cite{IKS1} applied to $\PP^2$.
Then,
we assume that $s\ge 2$ and
prove the inequality
\begin{equation}
N_{n(sF+E)}(\Sig_s)\ge n!\cdot
N_{n((s-1)F+E)}(\Sig_{s-1})\ ,\label{ne3}
\end{equation}
which
immediately
implies
$$\log N_{n(sF+E)}(\Sig_s)\ge\log
N_{nL}(\PP^2) +(s-1)\log n!\ ,$$
and hence
the inequality (\ref{ne2}) in view of
$\log N_{nL}(\PP^2)= 3n\log n+O(n)$.

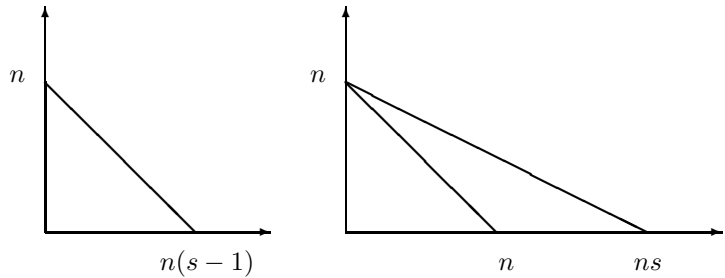
\begin{figure}\small
\setlength{\unitlength}{1cm}
\cl{\begin{picture}(11,5)(0,0)
\thinlines \put(1,1){\vector(1,0){3}}\put(1,1){\vector(0,1){3}}
\thicklines\put(1,1){\line(1,0){2}}\put(1,1){\line(0,1){2}}
\put(1,3){\line(1,-1){2}}\put(0.5,3){$n$}\put(2.5,0.5){$n(s-1)$}
\thinlines \put(5,1){\vector(1,0){5}}\put(5,1){\vector(0,1){3}}
\thicklines\put(7,1){\line(1,0){2}}\put(5,3){\line(1,-1){2}}
\put(5,3){\line(2,-1){4}}\put(4.5,3){$n$}\put(8.8,0.5){$ns$}\put(7,0.5){$n$}
\end{picture}}
\caption{Correspondence between rational curves on $\Sigma_{s-1}$
and $\Sigma_s$}\label{figgw1}
\end{figure}

To prove (\ref{ne3}),
notice that the number
$N_{n((s-1)F+E)}(\Sigma_{s-1})$
of rational curves in the linear system
$|n((s-1)F+E)|$
passing through
$(s+1)n-1$
generic
points in
$\Sigma_{s-1}$ can be viewed as the number of rational curves in
the linear system
$|n(sF+E)|$
on $\Sigma_s$ which pass
through
$(s+1)n-1$
generic points and have an ordinary
$n$--fold
singularity at some fixed point $z$
(this correspondence is provided by the birational transformation
$\Sigma_{s-1}\to\Sigma_s$ given in suitable affine
toric coordinates $x,y$ in
$\Sigma_{s-1}$ and $u,v$ in $\Sigma_s$ by $u=x, v=xy$; the
correspondence reflects an affine transformation of Newton
polygons as shown in Figure \ref{figgw1}). Choose now generic
points in a small neighborhood of $z$. The argument of the proof
of Lemma 5 in \cite{IKS1} confirms that any rational curve
$C\in|n(sF+E)|$ as above can be deformed inside the class of
rational curves passing through the initial $(s+1)n-1$ fixed
points, so that the $n$ local branches of $C$ at $z$ freely move
in transverse directions, and hence can be traced through the
newly chosen $n$ fixed points in an arbitrary order. Thus,
(\ref{ne3}) follows. \proofend

\begin{remark}
{\rm Lemma~\ref{nl1}
can also be proved using the tropical count. The
proof is completely similar to the proof of Theorem~3 (case of
$\PP^2$) in~\cite{IKS1}. One should just adapt the corresponding
lattice path constructed in~\cite{IKS1} to a triangle representing
the linear system $|n(sF+E)|$
on $\Sigma_s$}.
\end{remark}

\begin{lemma}\label{nl2-new}
Fix an integer $s \geq 1$.
Then, there exists an integer sequence $(T_n)$
which verifies the following properties:
$\log T_n=(s + 2)n\log n+O(n)$ and
for any $n$
there is a  generic collection $z_1, \ldots z_{(s+2)n-1}$
of $(s+2)n-1$
points in $\Sig_s$ such that
among the
rational curves
belonging to the
linear system
$|n(sF+E)|$
and  passing through $z_i$, $i=1,...,(s+2)n-1,$
at least $T_n$ curves
have only ordinary nodes as singular points and
intersect
each other
transversally outside of the points
$z_i$, $i=1,...,(s+2)n-1$.
\end{lemma}

{\bf Proof}\qua For $s=1$ the statement on ordinary nodes is
classical and holds for all the interpolating curves in the linear
system (see, for
example \cite{No} for a modern exposition). For $s>1$ one can
apply the construction described in the proof of Lemma \ref{nl1}
and observe that it preserves the statement on ordinary nodes. To
eliminate one by one eventual non-transversal intersections of
interpolating rational curves outside of the chosen points, it
suffices to move one of the chosen points along an interpolating
curve having a non-transversal intersection with another curve.
\proofend

\section{Proof of Theorem \ref{nt1} and Corollary \ref{cnt1}}

{\bf Proof of Theorem \ref{nt1}}\qua Let us first note that under the
hypotheses on~$D$ made in Theorem \ref{nt1}, the Gromov--Witten
invariants are enumerative. More precisely, deforming the complex
structure of $\PP^2_k$ to a generic almost complex one we observe,
first, that, due to $GW_D(\PP^2_k)\ne 0$ and the stability of the exceptional
divisors $E_1,\dots, E_k$ of $\PP^2_k\to\PP^2$,
the non-negativity of intersections, $D\cdot E_i\ge 0$,
holds,
and then, by Theorem 4.1 in
\cite{GP},
we obtain that due to $D\cdot c_1(\PP^2_k)\ >1$,
for a generic choice of the $k$ blown-up
points the number $\GW_{nD} (\PP^2_k), n\ge 1$, is equal to the
number $N_{nD} (\PP^2_k)$ of immersed irreducible rational curves
passing through $nD\cdot c_1(\PP^2_k) -1$ given generic points in
$\PP^2_k$.

The upper bound, $\log N_{nD}(\PP^2_k)\le
(D \cdot c_1(\PP^2_k))\, n\log n+O(n)$,
is given by Lemma 5 in
\cite{IKS1}.

To prove the lower bound, assume, first, that $D\cdot
c_1(\PP^2_k)>2$. By \cite{GP}, Theorem 4.1, we can represent $D$
by an immersion $\varphi\co \PP^1\to\PP^2_k$. We consider the bundle
${\cal N}_{\PP^1}=\varphi^*({\cal T}\PP^2_k)/{\cal T}\PP^1$ over
$\PP^1$ and compactify it into the rational ruled surface
$X=\Proj({\cal N}_{\PP^1}\oplus{\cal O}_{\PP^1})$. Pick K\"ahler
structures on $\PP^2_k$ and $X$ with the same periods on,
respectively, $D$ in $\PP^2_k$ and $\PP^1$ in $X$, and fix a
symplectic immersion~$\Phi$ of a small neighborhood ${\mathfrak
N}(\PP^1)\subset X$ into $\PP^2_k$ which extends $\varphi$. Such
an immersion~$\Phi$ exists due to the symplectic neighborhood
theorem (see, for example, \cite{Dusa}, Theorem 3.30) Notice that
$D\cdot c_1(\PP^2_k)=\PP^1\cdot c_1(X)$. Therefore,
applying~Lemma~\ref{nl2-new} to the linear system $|n\PP^1|$ we find in
it $T_n$, $T_n = D\cdot c_1(\PP^2_k)\, n\log n+O(n)$,
immersed rational curves which pass
through $D\cdot c_1(\PP^2_k)-1$ generic fixed points, have only
ordinary nodes as singular points, and intersect transversally
outside of the fixed points. Choose $0$ and $\infty$ in $\PP^1$ so
that the fibers over them are transversal to each of these
$T_n$
curves $C_i$ and do not contain any of the fixed
$D\cdot c_1(\PP^2_k)-1$ points. Now by vertical and horizontal
(toric) rescaling in $X$ we can make all the curves $C_i$ to be
$C^0$--close to $\PP^1$ and, moreover, $C^1$--close to it outside an
arbitrary small neighborhood of $0\in\PP^1$. As a consequence, we
get
$T_n$
immersed symplectic surfaces $\Phi(C_i)$ which
pass through some common $D\cdot c_1(\PP^2_k)-1$ points, have only
ordinary nodes as singular points, and are transversal to each
other outside of the common $D\cdot c_1(\PP^2_k)-1$ points.
Proceeding as in~\cite{Dusa2}, Lemma 3.2,
we construct a tamed almost complex
structure $J$ on $\PP^2_k$ for which all the surfaces $\Phi(C_i)$
are $J$--holomorphic (we start from neighborhoods of the common
points, where we retrieve a suitable almost complex structure from
$X$). Due to \cite{HLS}, the constructed $J$--holomorphic curves
represent discrete regular solutions of the interpolating problem.
Thus, to get the desired below bound it remains to notice that, as
it follows from \cite{HLS} and \cite{IS} (Corollaries 1.6 and
2.7), the space of generic almost complex structures is connected
and dense, and each regular solution counts for $+1$.

In the remaining case, $D\cdot c_1(\PP^2_k)=2$
and $D^2 > 0$, the conditions $D^2>0$
and $GW_D(\PP^2_k)\ne 0$ imply, by the standard gluing argument,
that $GW_{2D}(\PP^2_k)\ne 0$. Therefore, the preceding case
applies to $D'=2D$ and the lower bound now follows from the
monotonicity relation
\begin{equation}
N_{(n+1)D}(\PP^2_k)\ge N_{nD}(\PP^2_k)\ , n\ge 1.\label{ne4}
\end{equation}

To get (\ref{ne4}) we use again the
gluing of rational curves.
Namely, we construct an injective map from the set of rational
curves in $|nD|$ passing through $2n-1=nD \cdot c_1(\PP^2_k)-1$
fixed generic points to the set of rational curves in $|(n+1)D|$
passing through $2n+1$ generic points. Pick $2n$ generic points
$p_i$, $i=1,...,2n$, in $\PP^2_k$, and a rational curve
$C_1\in|D|$ passing through $p_{2n}$. We can assume that for any
curve~$C$ chosen among the rational curves belonging to $|nD|$ and
passing through $p_1,...,p_{2n-1}$, there is a point $z_C\in C\cap
C_1$ which is singular neither for $C$ no for $C_1$ and where the
curves~$C$ and~$C_1$ intersect transversally. Pick a generic point
$p'\in C_1$ and a point $p_{2n+1}\not\in C_1$ in a small
neighborhood of $p'$. Then, there exists a one-parameter
deformation of $C\cup C_1$ consisting of rational curves in
$|(n+1)D|$ and such that the point $z_C$ smoothes out and the
points $p_1,...,p_{2n}$ remain fixed (see, for example, \cite{GK},
Proposition 5.2, or \cite{Ko}, Ch.~II, Theorem 7.6). This family
sweeps a neighborhood of $p'$, and hence we obtain a rational
curve $C'\in|(n+1)D|$ passing through
$p_1,...,p_{2n-1},p_{2n},p_{2n+1}$. \proofend

\begin{remark}
{\rm
The relation~(\ref{ne1}) of the statement of
Theorem~\ref{nt1} is also valid for a homology class $D\in
H_2(\PP^2_k;\Z)$ such that $GW_D(\PP^2_k) > 1$, $D \cdot
c_1(\PP^2_k) = 1$, and $D^2 > 0$.
Indeed,
putting $D' = 2D$, we have $D'\cdot c_1(\PP^2_k) = 2$ and $(D')^2
> 0$. By \cite{GP}, Theorem 4.1, we can find in the linear system
$|D|$ two distinct immersed rational curves. Then, deforming these
curves as in the proof of Theorem~\ref{nt1}, we get
$GW_{D'}(\PP^2_k) > 0$. Thus, the relation~(\ref{ne1}) holds for $D'$,
and we can prove this relation for~$D$ using
the
inequality (\ref{ne4}) as is done in the proof of
Theorem~\ref{nt1}.}
\end{remark}

{\bf Proof of Corollary \ref{cnt1}}\qua The first
statement immediately follows from Theorem 1 and the equality
$\GW_{nD} (\PP^2_k)=N_{nD} (\PP^2_k)$ explained in the beginning
of the proof of Theorem 1.

To
prove the second statement, we observe first
that any ample divisor~$D$
on $\PP^2_k$, $k\le 9$, is represented by a nodal rational curve.
For $k=1$ or $2$ this is trivial. For $3\le k\le 9$ this
follows from
Theorem 5.2 in
\cite{GLS}, which states the existence of rational
nodal curves in certain linear systems in $\PP^2_k$. Indeed,
given an expansion $D = dL-d_1E_1-...-d_kE_k$
for any base $(L,E_1,...,E_k)$ of
$\Pic(\PP^2_k)$
satisfying $L^2=-E_1^2=...=-E_k^2=1$, the
ampleness of~$D$
yields that $d,d_1,...,d_k>0$. Furthermore, by base
changes in $\Pic(\PP^2_k)$ induced by Cremona transformations
(see \cite{GLS}, section 5.1), we can achieve
$d\ge\max_{i\ne j\ne l}(d_i+d_j+d_l)$, which is the minimality
condition of Theorem 5.2 in~\cite{GLS}.
At last, the remaining
condition of this
theorem, $3d>d_1+...+d_k$ (which can be also written
as
$D \cdot c_1(\PP^2_k) > 0$),
follows from the positivity of
intersection with the strict transform of a plane cubic
passing
through the blown-up points.

Since the existence of a nodal rational curve in
$\vert D\vert$ implies $GW_D(\PP^2_k)\ne 0$,
and
since, in addition, $D^2 > 0$ for any ample $D$, the second
statement of the corollary is proved for ample divisors $D$
satisfying
$D \cdot c_1(\PP^2_k) \ge 2$.

Now
let us consider the case
$D \cdot c_1(\PP^2_k) = 1$,
and put $D' = 2D$. We have \mbox{$D' \cdot c_1(\PP^2_k) = 2$} and
$(D')^2
> 0$, and once more by Theorem 5.2 in~\cite{GLS} we get
$GW_{D'}(\PP^2_k)\ne 0$. Hence, the relation~(\ref{ne10}) holds for
$D'$, and finally we deduce this relation for~$D$ using
the inequality (\ref{ne4}) as is done in the proof of
Theorem~\ref{nt1}. \proofend

\section{Welschinger invariants
of real rational surfaces}\label{last}

Recall that the
Welschinger invariants depend not only on a
homology class,
but also on a number of
non-real
points in a real
configuration of points (see \cite{W,W1} for the definition and
properties of the Welschinger invariants).
Denote
by
$W_{nD}(\PP^2_k)$ the Welschinger invariant
which counts, with weights $\pm 1$, the real rational curves
belonging to the linear system $|nD|$ and passing through $nD\cdot
c_1(\PP^2_k)-1$ given generic real points in $\PP^2_k$. As is proved in
\cite{IKS1}, in the cases $k=1,2,3$, the
same
relation
as (\ref{ne1}) holds for $W_{nD}(\PP^2_k)$
(instead of $\GW_{nD}
(\PP^2_k)$).
This motivates the following conjecture.

\begin{conjecture}\label{nc1}
Assume that $\PP^2_k$ is obtained from $\PP^2$ by blowing up $k$
generic real points
and is equipped with its
natural real structure.
Let $D\subset
\PP^2_k$ be a real ample divisor. Then, the Welschinger invariants
$W_{nD}
(\PP^2_k)$
satisfy the
relation
$$\lim_{n\to
+\infty} \frac{\log W_{nD}(\PP^2_k)}{n\log n}=D\cdot c_1({\PP^2_k})\ .$$
\end{conjecture}

One could try to prove Conjecture~\ref{nc1}
using the same construction as
in the proof of Theorem~\ref{nt1}.
However, this approach does not give immediately the result,
since a real regular solution
to the interpolation problem contributes $\pm 1$
to the Welschinger invariant.
Thus, to get an asymptotic lower bound,
it is not enough to present an appropriate number
of interpolating real rational curves.

\subsection*{Acknowledgements}

The idea of this work came to the authors during their stay at the
Max-Planck-Institut f\"ur Mathematik in Bonn, and we thank this
institution for the hospitality and excellent work conditions. The
first and the second authors are members of Research Training Network
RAAG CT-2001-00271.  The third author was supported by the Israel
Science Foundation grant no. 465/04 and by the
Nermann--Minkowski--Minerva Center for Geometry at the Tel Aviv
University.

We are also grateful to the referee for useful
remarks.

\end{document}

%% file: gtoutput.tex
\def\ifplaintex{\expandafter\ifx\csname documentclass\endcsname\relax}

\ifplaintex 
\else
\fi

\expandafter\ifx\csname beginpicture\endcsname\relax
\expandafter\ifx\csname documentclass\endcsname\relax
\input pictex \else
\input prepictex \input pictex \input postpictex \fi\fi

\def\gt{{\mathsurround=0pt\it $\cal G\mskip-2mu$eometry \&\ 
$\cal T\!\!$opology}}        

\def\gtp{{\mathsurround=0pt\it $\cal G\mskip-2mu$eometry \&\ 
$\cal T\!\!$opology $\cal P\!$ublications}}  

\def\lognumber#1{\def\thelognumber{#1}}
\def\volumenumber#1{\def\thevolumenumber{#1}}
\def\papernumber#1{\def\thepapernumber{#1}}
\def\volumeyear#1{\def\thevolumeyear{#1}}

\def\pagenumbers#1#2{\def\startpage{#1}\def\finishpage{#2}}
\def\published#1{\def\publishdate{#1}}
\def\proposed#1{\def\theproposer{#1}}
\def\seconded#1{\def\theseconders{#1}}
\def\received#1{\def\receiveddate{#1}}

\def\accepted#1{\def\accepteddate{#1}}
\def\asciititle#1{\def\theasciititle{#1}}

\def\asciiauthors#1{\def\theasciiauthors{#1}}
\def\asciiaddress#1{\def\theasciiaddress{#1}}
\def\asciiemail#1{\def\theasciiemail{#1}}

\long\def\asciiabstract#1{\long\def\theasciiabstract{#1}}
\def\asciikeywords#1{\def\theasciikeywords{#1}}
\def\shortauthors#1{\def\theshortauthors{#1}}

\let\\\par\let\thelognumber\relax
\let\thevolumenumber\relax\let\thepapernumber\relax
\let\thevolumeyear\relax\let\thesamplenumber\relax\let\startpage\relax
\let\finishpage\relax\let\publishdate\relax\let\receiveddate\relax
\let\reviseddate\relax\let\accepteddate\relax\let\theasciititle\relax
\let\theasciiauthors\relax\let\theasciiaddress\relax
\let\theasciiabstract\relax\let\theasciikeywords\relax
\let\theasciiemail\relax\let\theshortauthors\relax\let\theshorttitle\relax

\long\def\maketitlep{   

\count0=\startpage

\gt\hfill      
\beginpicture
\setcoordinatesystem units <0.33truein, 0.33truein> point at 2.2 0.9
\setplotsymbol ({$\cal G$})
\plotsymbolspacing=9truept
\circulararc 315 degrees from 0 1 center at 0 0
\setplotsymbol ({$\cal T$})
\circulararc 315 degrees from 1 -1 center at 1 0
\endpicture
\break
{\small\ifx\thesamplenumber\relax 
Volume \else Sample
\fi\thevolumenumber\ (\thevolumeyear)
\startpage--\finishpage\nl
Published: \publishdate}
\vglue 0.5truein plus 0.4fil minus 0.1truein

{\parskip=0pt\leftskip 0pt plus 1fil\def\\{\par\smallskip}{\ifplaintex\large
\else\Large\fi\bf\thetitle}\par\medskip}   

\vglue 0pt plus 0.1fil 

{\parskip=0pt\leftskip 0pt plus 1fil\def\\{\par}{\sc\theauthors}
\par\medskip}

\vglue 0pt plus 0.1fil 

{\small\parskip=0pt\let\newline\\
{\leftskip 0pt plus 1fil\def\\{\par}{\sl\theaddress}\par}
\expandafter\ifx\theemail\relax    
\relax\else\vglue 5pt plus 0.02fil minus 2pt\def\\{\stdspace{\rm 
and}\stdspace} 
\cl{Email:\stdspace\tt\theemail}\fi
\ifx\theurl\relax                  
\relax\else\vglue 5pt plus 0.02fil minus 2pt\def\\{\stdspace{\rm 
and}\stdspace}
\cl{URL:\stdspace\tt\theurl}\fi\par}

\vglue 7pt plus 0.3fil minus 3pt

{\bf Abstract}
\vglue 5pt plus 0.1fil minus 2pt

\theabstract

\vglue 7pt plus 0.3fil minus 3pt

{\bf AMS Classification numbers}\quad Primary:\quad \theprimaryclass

Secondary:\quad \thesecondaryclass

\vglue 5pt plus 0.3fil minus 2pt

{\bf Keywords:}\quad \thekeywords

\vglue 10pt plus 0.5fil minus 5pt

{\small  Proposed: \theproposer\hfill Received: \receiveddate\nl
Seconded: \theseconders\hfill 
\ifx\reviseddate\relax                         
Accepted: \accepteddate                        
\else
Revised: \reviseddate                          
\fi}
\eject
}       

\let\maketitlepage\maketitlep
\let\maketitle\maketitlepage

\font\phead=cmsl9 scaled 950
\font=cmsl9 scaled 1050
\font\pnum=cmbx10 scaled 913
\font\lnum=cmbx10 
\font\pfoot=cmsl9 scaled 950
\font=cmsl9 scaled 1050
\ifplaintex
\headline{\vbox to 0pt{\vskip -4.5mm\line{\small\phead\ifnum
\count0=\startpage ISSN 1364-0380 (on line)
1465-3060 (printed) \hfill {\pnum\folio}\else\ifodd\count0\def\\{ }%
\ifx\theshorttitle\relax\thetitle\else\theshorttitle\fi\hfill{\pnum\folio}
\else\def\\{ and }{\pnum\folio}\hfill\ifx\theshortauthors\relax\theauthors
\else\theshortauthors\fi\fi\fi}\vss}}
\footline{\vbox to 0pt{\vglue 0mm\line{\small\pfoot\ifnum\count0=\startpage
\copyright\ \gtp\hfill\else
\gt, Volume \thevolumenumber\ (\thevolumeyear)\hfill\fi}\vss
}}
\else
\makeatletter
\def\@oddhead{{\small\lhead\ifnum\count0=\startpage ISSN 1364-0380 (on line)
1465-3060 (printed) \hfill {\lnum\number\count0}\else\ifodd\count0
\def\\{ }\ifx\theshorttitle\relax \thetitle \else\theshorttitle\fi\hfill
{\lnum\number\count0}\else\def\\{ and }{\lnum\number\count0}
\hfill\ifx\theshortauthors\relax 
\theauthors\else\theshortauthors\fi\fi\fi}}\def\@evenhead{\@oddhead}
\def\@oddfoot{\small\lfoot\ifnum\count0=\startpage\copyright\ \gtp\hfill\else
\gt, Volume \thevolumenumber\ (\thevolumeyear)\hfill\fi}
\def\@evenfoot{\@oddfoot}
\makeatother
\fi

\newwrite\gtoutfile
\long\gdef\makeheadfile{  
{\def\\{, }\def\s{ }
\immediate\openout\gtoutfile head.xxx
\immediate\write\gtoutfile{Proxy-for: \ifx\theasciiauthors\relax
\theauthors\else\theasciiauthors\fi\s<\ifx\theasciiemail\relax\theemail\else\theasciiemail\fi>}
\immediate\write\gtoutfile{\noexpand\\}
\immediate\write\gtoutfile{Authors: \ifx\theasciiauthors\relax
\theauthors\else\theasciiauthors\fi}
{\def\\{ }\immediate\write\gtoutfile{Title: \ifx\theasciititle\relax
\thetitle\else\theasciititle\fi}}
\immediate\write\gtoutfile{Subj-class: GT or SG or MG etc}
\immediate\write\gtoutfile{MSC-class: \theprimaryclass\ifx\thesecondaryclass\relax\else, \thesecondaryclass\fi}
\immediate\write\gtoutfile{Journal-ref: Geom. Topol. \thevolumenumber
(\thevolumeyear) \startpage-\finishpage}
\immediate\write\gtoutfile{Comments: Published by Geometry and Topology at}
\immediate\write\gtoutfile{\s\s http://www.maths.warwick.ac.uk/gt/GTVol\thevolumenumber/paper\thepapernumber.abs.html}
\immediate\write\gtoutfile{\noexpand\\}
\immediate\write\gtoutfile{}
\ifx\theasciiabstract\relax
\immediate\write\gtoutfile{\theabstract}\else
\immediate\write\gtoutfile{\theasciiabstract}\fi
\immediate\write\gtoutfile{}
\immediate\write\gtoutfile{\noexpand\\}
\immediate\write\gtoutfile{}
\immediate\closeout\gtoutfile}}  

\def\maketitlepage{\maketitlep\makeheadfile}
\let\maketitle\maketitlepage